\documentclass[11pt]{amsart}
\usepackage{amsmath,amsthm,amssymb}
\newtheorem*{definition}{Definition}
\newtheorem{theorem}{Theorem}
\newtheorem*{lemma}{The right-multiplication lemma}
\newtheorem*{mercer}{Mercer's theorem}
\newtheorem{remark}{Remark}
\newcommand{\Int}{\int\limits}
\newcommand{\Sum}{\sum\limits}
\newcommand{\Sup}{\sup\limits}

\DeclareMathOperator{\RE}{Re}
\DeclareMathOperator{\IM}{Im}
\def\R{\mathbb R_+}
\def\bsy{\boldsymbol}\def\le{\leqslant}\def\ge{\geqslant}
\def\adj#1{\overline{#1}}
\begin{document}

\title{An extension of Mercer's theorem to unbounded operators}

\author{I. M. Novitski\u i}
\address{Institute for Applied Mathematics, Russian Academy of Sciences,
9, Shevchenko Street, Khabarovsk 680 000, Russia}
\email{novim@iam.khv.ru}

\author{M. A. Romanov}
\address{Institute for Applied Mathematics, Russian Academy of Sciences,
9, Shevchenko Street, Khabarovsk 680 000, Russia}
\email{romanov@iam.khv.ru}

\date{December 31, 1997.}

\begin{abstract}
We give some extensions of Mercer's theorem to continuous Carleman kernels
inducing unbounded integral operators.
\end{abstract}

\maketitle

\section*{Introduction}

The source of the following theorem is \cite{12}.
\begin{mercer}
Let $T$ be a positive, integral operator on $L_2[a,b]$ with continuous kernel
$K(s,t)=\overline{K(t,s)}$ on $[a,b]^2$ $(|a|,|b|<\infty)$. Then the kernel
$K(s,t)$ can be represented by the bilinear series
$$
K(s,t)=\Sum_{n=1}^\infty\lambda_n\varphi_n(s)
\overline{\varphi_n(t)}
$$
absolutely and uniformly convergent on $[a,b]^2$, where $\lambda_n\ge0$
$(n=1,2,3,\dots)$ are the eigenvalues of operator $T$ and $\varphi_n$ 
$(n=1,2,3,\dots)$ are the corresponding orthonormal eigenfunctions.
\end{mercer}

Within modern theory of integral operators, the theorem just formulated is
related to the problem of representing the kernel via spectrum and
eigenfunctions (see, for example, \cite{9}). In the works \cite{2},\ \cite{4},\ \cite{5},\ 
\cite{11},\ \cite{13}-\cite{15}, Mercer's theorem has been generalized to wider 
classes of kernels inducing bounded integral operators.
In the present paper, we study the problem of spectral representing
for a continuous kernel in that case when the kernel induces an unbounded
normal operator whose spectrum is in a sector of angle less than $\pi$ with
vertex at $0$.

The following definitions and notions are needed in what follows.
Let $\R=[0,\infty)$, and let $L_2$ be the Hilbert space of complex-valued
measurable functions on $\R$ that are square integrable in the Lebesgue sense,
equipped with the inner product
$$
\langle f,g\rangle=\int_{\R} f(s)\overline{g(s)}\,ds
$$
and the norm $\Vert f\Vert=\langle f,f\rangle^{\frac{1}2}$.

A linear operator $T:D_T\rightarrow L_2$, where $D_T$ is a dense linear manifold
in $L_2$, is said to be {\it integral\/} if there exists a
measurable almost-everywhere finite-valued function $K(s,t)$ on
$\R^2$, a {\it kernel}, such that, for every $f\in D_T$,
$$(Tf)(s)=\int_{\R} K(s,t)f(t)\,dt$$
for almost every $s\in\R$.

A kernel $K(s,t)$ is said to be {\it Carleman} if $K(s,\cdot)\in L_2$
for almost every fixed $s\in\R$.

An integral operator induced by a Carleman kernel is called {\it Carleman operator},
and it is called {\it bi-Carleman} if $K(\cdot,t)\in L_2$ for almost every $t\in\R$.
In the last case the kernel induces two {\it Carleman functions\/}
from $\R$ to $L_2$ by $k(s)=\overline{K(s,\cdot)}$,
$k^\ast(t)=K(\cdot,t)$ for all those $s$, $t\in\R$
for which $K(s,\cdot)\in L_2$ and $K(\cdot,t)\in L_2$.
We need the following important result of Carleman operator theory

\begin{lemma}
{\rm (\cite{8}).} Let $T:D_T\to L_2$ be a Carleman operator, and let $R:L_2\to L_2$ be
a bounded operator. Then the product $TR$ is a Carleman operator with Carleman function
$R^*(k(s))$, where $k(s)$ is Carleman function of $T$.
\end{lemma}

If $X$ is a locally compact space, and if $B$ is a Banach space with norm
$\Vert \cdot\Vert_B$, then let $C(X,B)$ denote the Banach space (with the norm
$\Vert f\Vert_{C(X,B)}=\Sup_X\Vert f(x)\Vert_B$) of all continuous functions
from $X$ to $B$ vanishing at infinity (the latter means that for each
fixed $f\in C(X,B)$, given an arbitrary $\varepsilon>0$, there exists a compact set
$X(\varepsilon,f)\subset X$ such that $\Vert f(x)\Vert_B<\varepsilon$
if $x\not\in X(\varepsilon,f)$).

\begin{definition}
{\rm (\cite{15,16}).} A kernel $K$ of a bi-Carleman operator $T$ is called a $K^0$-kernel, if
$K\in C(\R^2,\mathbb C)$ and its Carleman functions $k$, $k^\ast\in C(\R,L_2)$.
\end{definition}

In Section 1, we give an ``integral'' analogue of Mercer's theorem (Theorem~1) that represents
``functions of a $K^0$-kernel'' by uniformly convergent principal value
Lebesgue-Stieltjes integrals with respect to the spectral function of this kernel.
In Section~2, $K^0$-kernels of unbounded diagonal operators are considered.
There we prove Theorem~2 on representations of $K^0$-kernels by absolutely and
uniformly convergent bilinear series. In the case of bounded integral operators
Theorems~1 and 2 were proved in \cite{14,15}.

\section{Integral representations of $K^0$-kernels}

Let $N:D_N\to L_2$ be a normal operator, that is, a closed linear operator
with a dense in $L_2$ domain $D_N = D_{N^*}$ and such that $NN^* = N^*N$.
Let $E(\cdot)$ be the resolution of identity induced by $N$ (see \cite{18}).
Assume that $N$ is an integral operator having a $K^0$-kernel
$N(s,t)$ and Carleman functions $\nu(s)=\overline{N(s,\cdot)}$,
$\nu^*(t)=N(\cdot,t)$.

Let $\Omega$ be the $\sigma$-algebra of Borel sets in $\mathbb C$, and let $\Omega_0$
be the set of all $\omega\in\Omega$ whose closures in $\mathbb C$ do not
contain $0$. If $\chi_\omega$ is the characteristic function of the set
$\omega\in\Omega_0$, then $\chi_\omega (z)=zv_\omega (z)$, where
$v_\omega (z)=\chi_\omega (z)/z$. From the multiplicative property of the spectral
measure $E(\cdot)$ it follows that $E(\omega)=\chi_\omega (N)=Nv_\omega (N)$
(here and throughout, $f(N)$, given an $\Omega$-measurable function $f$,
denote a normal operator defined by
$$
\langle f(N)h,g\rangle=\int_\mathbb C f(z)\langle E(dz)h,g\rangle\quad
\text{for all $h$, $g\in D_{f(N)}$}
$$
with domain
$$
D_{f(N)}=\left\{
x\in L_2:\int_\mathbb C |f(z)|^2\langle E(dz)x,x\rangle<\infty
\right\}).
$$
Since $B(\omega)=v_\omega (N)$ is a bounded operator for each
$\omega\in\Omega_0$, it follows by the right-multiplication lemma that the
orthogonal projection $E(\omega)$ is a bounded bi-Carleman operator,
with the Carleman functions
$$
e(s;\omega)=(B(\omega))^*(\nu(s)),\quad e^*(t;\omega)=e(t;\omega)\eqno(1.0)
$$
belonging to $C(\R,L_2)$. Since $E^2(\omega)=E(\omega)$, the kernel
of $E(\omega)$ can be computed as follows
$$
E(s,t;\omega)=\langle e(t;\omega),e(s;\omega)\rangle
$$
for all $(s,t)\in\R^2$, and hence belongs to $C(\R^2,\mathbb C)$.

\begin{definition}
A set function $E(s,t;\cdot):\Omega_0\rightarrow C(\R^2,\mathbb C)$ whose value
on every $\omega\in\Omega_0$ is the inducing $K^0$-kernel $E(s,t;\omega)$
of $E(\omega)$ is called a {\it spectral function\/} for the kernel $N(s,t)$.
\end{definition}

\begin{remark}
For other definitions of spectral functions of Carleman kernels we refer to
the works \cite{6,3,10,19,1,8}.
\end{remark}

Let $\omega_\varepsilon=\mathbb C\setminus\{z\in\mathbb C:\,\vert
z\vert\le\varepsilon,\,\varepsilon>0\}$, and let $\bsy\Phi$ denote the
family of all functions $\bsy\phi$ of the form $\bsy\phi(z)=zv(z)$ where $v(z)$ is
a bounded $\Omega$-measurable function on $\mathbb C$.
Let $\Omega_\varepsilon$ be the algebra of all Borel subsets of $\omega_\varepsilon$.
Since $E(\omega_\varepsilon)E(\sigma)=E(\sigma)$
for each $\sigma\in\Omega_\varepsilon$, the right-multiplication lemma yields
$e(s;\sigma)=E(\sigma)(e(s;\omega_\varepsilon))$ and hence
$$
E(s,t;\sigma)=\langle e(t;\sigma),e(s;\sigma)\rangle=
\langle E(\sigma)(e(t;\omega_\varepsilon)),e(s;\omega_\varepsilon)\rangle
\eqno(1.1)
$$
for all $s$, $t\in\R$.
Whence it follows that
$$
\Vert E(\cdot,\cdot;\sigma)\Vert_{C(\R^2,\mathbb C)}\leq
\|e(\cdot;\omega_\varepsilon)\|_{C(\R,L_2)}^2
$$
for each $\sigma\in\Omega_\varepsilon$.
The last inequality implies that, for each
$\varepsilon>0$, the spectral function $E(s,t;\cdot)$ is bounded on
$\Omega_\varepsilon$. Furthermore, it is additive; this property follows from
that of $E(\cdot)$. Therefore, for each $\bsy\phi\in\bsy\Phi$,
the Lebesgue-Stieltjes integrals
$$
\varPhi_\varepsilon(s,t)=\int_{\omega_\varepsilon}\bsy\phi(z)E(s,t;dz),
\eqno(1.2)
$$
can be formed for every $\varepsilon>0$ and for all $s$, $t\in\R$.

\begin{theorem}
Let $N:D_N\to L_2$ be a normal, integral operator induced by a $K^0$-kernel
$N(s,t)$, the spectrum of which is in a sector $\mathfrak S$ of angle less than $\pi$
with vertex at $0$. Then for every $\bsy\phi\in\bsy\Phi$ operator
$\bsy\phi(N)$ is also integral operator with $K^0$-kernel. Moreover, the
$K^0$-kernel $\varPhi(s,t)$ inducing $\bsy\phi(N)$ can be represented by the
principal value integral with singularity at $z=0$:
$$
\varPhi(s,t)=\int_\mathbb C\bsy\phi(z)\,E(s,t;\,dz)=\lim_{\varepsilon\to 0}
\int_{\omega_{\varepsilon}}\bsy\phi(z)\,E(s,t;\,dz)
$$
for all $s$, $t\in\R$, where the integral converges to $\varPhi(s,t)$
in $C(\R^2,\mathbb C)$ as $\varepsilon\rightarrow 0$ along arbitrary decreasing
sequence of positive numbers.
\end{theorem}
\begin{proof}
From multiplicative property it follows that
$$
\langle\bsy\phi(N)f,h\rangle=\Int_\mathbb C zv(z)\langle E(dz)f,h\rangle=
\langle Nv(N)f,h\rangle\quad\text{for all}\ f,h\in D_{\bsy\phi(N)},
$$
that is, $\bsy\phi(N)=Nv(N)$ is a Carleman operator with a kernel $\varPhi(s,t)$
and a Carleman function
$$
\varphi(s)=\overline{\varPhi(s,\cdot)}=(v(N))^*(\nu(s)),\eqno(1.3)
$$
which is in $C(\R,L_2)$.
Since $\bsy\phi(N)$ is a normal operator, it is bi-Carleman
(see Corollary 2.19 from \cite[p. 131]{8}). In addition, the Carleman function of
the adjoint
$$
\langle(\bsy\phi(N))^*f,h\rangle=
\Int_\mathbb C\overline z\overline{v(z)}\langle E(dz)f,h\rangle=
\langle N^*(v(N))^*f,h\rangle,\quad f,h\in D_{\bsy\phi(N)}
$$
has, by the right-multiplication lemma, the form
$$
\varphi^*(t)=\varPhi(\cdot,t)=v(N)(\nu^*(t))\quad(t\in\R)
$$
and hence belongs to $C(\R,L_2)$ too.
Using (1.0) and (1.1), compute the integral (1.2) as follows
$$
\gathered
\varPhi_\varepsilon (s,t)=\int_{\omega_\varepsilon}\bsy\phi (z)\langle
E(dz)(e(t;\omega_\varepsilon)),e(s;\omega_\varepsilon)\rangle=\\
=\int_{\omega_\varepsilon}\bsy\phi(z)\langle
E(dz)(B(\omega_\varepsilon))^*(\nu(t)),
(B(\omega_\varepsilon))^*(\nu(s))\rangle=\\
=\langle\bsy\phi(N)B(\omega_\varepsilon)(B(\omega_\varepsilon))^*
(\nu(t)),\nu(s)\rangle
\endgathered\eqno(1.4)
$$
for all $s$, $t\in\R$. It is clear that the function $\varPhi_\varepsilon(s,t)$ 
belongs to $C(\R^2,\mathbb C)$ and that this is a kernel of the integral operator
$\bsy\phi(N)E(\omega_\varepsilon)=NE(\omega_\varepsilon)v(N)$ with Carleman
function
$$
\varphi_\varepsilon(s)=\overline{\varPhi_\varepsilon (s,\cdot)}=
(v(N))^*E(\omega_\varepsilon)(\nu(s))\quad(s\in\R)\eqno(1.5)
$$
belonging to $C(\R,L_2)$.
The function $E(\{0\})(\nu(s))$ from $C(\R,L_2)$ is equal identically to zero
on $\R$, because it is a Carleman function of the integral operator $NE(\{0\})$,
which, by the multiplicative property, is the null operator:
$$
\langle NE(\{0\})f,h\rangle=\Int_\mathbb C z\chi_{\{0\}}(z)
\langle E(dz)f,h\rangle=0,\quad f,h\in L_2.
$$
Let $\{\varepsilon_n\}\subset\R$ be an arbitrary sequence
decreasing to $0$, and let $\omega_0=\mathbb C\setminus\{0\}$.
By virtue of (1.3) and (1.5),
$$
\gathered
\lim_{n\to\infty}
\|\varphi-\varphi_{\varepsilon_n}\|_{C(\R,L_2)}\le\\
\le\|(v(N))^*\|\,\lim_{n\to\infty}\|
(E(\omega_0)-E(\omega_{\varepsilon_n}))\nu\|_{C(\R,L_2)}=0,
\endgathered\eqno(1.6)
$$
because $E(\omega_{\varepsilon_n})\to E(\omega_0)$ whenever $n\to\infty$
in strong operator topology and the set $\left\{\nu(s):s\in\R\right\}$ is precompact
in $L_2$ (see \cite[p. 193]{7}). The following property
$$\lim_{n\to\infty}
\|\varphi^*-\varphi_{\varepsilon_n}^*\|_{C(\R,L_2)}=0,\eqno(1.6')$$
where $\varphi_{\varepsilon_n}^*(t)=\varPhi_{\varepsilon_n}(\cdot,t)$, 
$n=1,2,3,\dots$ can be proved analogously.

Without loss of generality, assume that the sector $\mathfrak S$ is bounded by
the rays $\IM z=\pm l\RE z$, $\RE z\ge 0$, $l>0$.
It is obvious that
$$
\Sup_{z\in \mathfrak S}|v_1(z)|\le l,\ \text{where}\
v_1(z)=\dfrac{\IM z}{\RE z}.\eqno(1.7)
$$
Consider a self-adjoint integral operator $X=(N+N^*)/2$ induced by the $K^0$-kernel
$X(s,t)=\overline{X(t,s)}=\big(N(s,t)+\overline{N(t,s)}\big)/2$ having the
Carleman function $x(s)=\overline{X(s,\cdot)}$.
Fix $\varepsilon>0$ and consider the bounded integral operator
$X\big(I-E(\omega_\varepsilon)\big)$ with a $K^0$-kernel
$X(s,t)-X_\varepsilon(s,t)$, where
$$
\gathered
X_\varepsilon (s,t)=\int_{\mathfrak S\cap\omega_\varepsilon}
\RE z E(s,t;dz)=\\
=\langle XB(\omega_\varepsilon)(B(\omega_\varepsilon))^*
(\nu(t)),\nu(s)\rangle
\quad\left((s,t)\in\R^2\right)
\endgathered\eqno(1.8)
$$
is the $K^0$-kernel of integral operator $XE(\omega_\varepsilon)$ inducing the
Carleman function $x_\varepsilon(s)=\overline{X_\varepsilon(s,\cdot)}$.
This operator is positive, since
$$
\langle X(I-E(\omega_\varepsilon))f,f\rangle=
\Int_{\mathfrak S}\RE z\chi_{\mathbb C\setminus\omega_\varepsilon}(z)
\langle E(dz)f,f\rangle\ge 0,\quad f\in L_2.
$$
Hence, its $K^0$-kernel satisfies the inequality
$X(s,s)-X_\varepsilon (s,s)\ge 0$ for all $s\in\R$
(see\ \cite[p. 98]{17}), whence
$$
X_\varepsilon (s,s)\le X(s,s)\quad\text{for all $s\in\R$}.\eqno(1.9)
$$
The bounded self-adjoint operator
$$
x(\varepsilon_m,\varepsilon_n)=
X\left(B(\omega_{\varepsilon_m})
(B(\omega_{\varepsilon_m}))^*
-
B(\omega_{\varepsilon_n})
(B(\omega_{\varepsilon_n}))^*\right)
$$
is positive if $0<\varepsilon_m\le\varepsilon_n$:
$$
\langle  x(\varepsilon_m,\varepsilon_n)f,f\rangle=
\Int_{{\mathfrak S}\cap(\omega_{\varepsilon_m}\setminus
\omega_{\varepsilon_n})}\frac{\RE z}{|z|^2}\langle E(dz)f,f\rangle\ge 0,
\quad f\in L_2.
$$
Substituting $f=\nu(s)$ to the last inequality and taking into account
(1.8) we obtain
$X_{\varepsilon_m}(s,s)\ge X_{\varepsilon_n}(s,s)$ for all $s\in\R$.
Conclude via (1.9) that the sequence $X_{\varepsilon_n}(s,s)$ ($n=1,2,3,\dots$)
of functions from $C(\R,\mathbb C)$ converges in $\R$ if
$\varepsilon_n\searrow 0$ as $n\to\infty$.
Apply the generalized Schwarz inequality \cite[p. 78]{1} to the 
positive operator $x(\varepsilon_m,\varepsilon_n)$ to write
$$
\gathered
|X_{\varepsilon_m}(s,t)-X_{\varepsilon_n}(s,t)|^2=\\
=|\langle x(\varepsilon_m,\varepsilon_n)(\nu(t)),\nu(s)\rangle|^2
\le\langle x(\varepsilon_m,\varepsilon_n)(\nu(s)),\nu(s)\rangle
\langle x(\varepsilon_m,\varepsilon_n)(\nu(t)),\nu(t)\rangle=\\
=\left( X_{\varepsilon_m}(s,s)-X_{\varepsilon_n}(s,s)\right)
\left( X_{\varepsilon_m}(t,t)-X_{\varepsilon_n}(t,t)\right).
\endgathered
$$
Whence, by (1.9), one can infer that
$$
\lim_{n\to\infty}\Sup_{t\in\R}
|X_{\varepsilon_n}(s,t)-X'(s,t)|=0\quad
\text{for all\ } s\in\R,
$$
where $X'(s,\cdot)\in C(\R,\mathbb C)$ for each $s\in\R$.
Now apply the preceding arguments to the functions
$x$, $x_{\varepsilon_n}$, $\dfrac{\RE z}{z}$ instead of
$\varphi$, $\varphi_{\varepsilon_n}$, $v(z)$
respectively to conclude that
$$
\lim_{n\to\infty}\|x_{\varepsilon_n}-x\|_{C(\R,L_2)}=0.
$$
The functions $X(s,t)$ and $X'(s,t)$ are continuous with respect to $t$  
and coincide for each fixed $s\in\R$ and all $t\in\R$, because from two
last limit relations it follows that, for every $s\in\R$,
$$
X(s,t)=X'(s,t)
$$
for almost every $t\in\R$
(see \cite[p. 42]{1}). Assume that $X(\infty,\infty)=0$ and that
$X_{\varepsilon_n}(\infty,\infty)=0$ ($n=1,2,3,\dots\,$), and
conclude, by the Dini theorem, that the monotone increasing sequence
of $\{X_{\varepsilon_n}(s,s)\}$ uniformly converges on the compactum
$[0,\infty]$ to the function $X(s,s)$. Consequently, one can write
$$
\lim_{n\to\infty}\Sup_{s\in\R}
|X_{\varepsilon_n}(s,s)-X(s,s)|=0.\eqno(1.10)
$$
Consider the bounded operator
$$
\bsy\phi(\varepsilon_m,\varepsilon_n)=
\bsy\phi(N)
\left(B(\omega_{\varepsilon_m})
(B(\omega_{\varepsilon_m}))^*
-
B(\omega_{\varepsilon_n})
(B(\omega_{\varepsilon_n}))^*\right),
\quad
0< \varepsilon_m\le\varepsilon_n.
$$
By (1.7), we have
$$
\gathered
\langle\bsy\phi(\varepsilon_m,\varepsilon_n)f,h\rangle=
\Int_{\mathfrak S\cap(\omega_{\varepsilon_m}\setminus\omega_{\varepsilon_n})}
v(z)\big(1+i\cdot v_1(z)\big)
\frac{\RE z}{|z|^2}\langle E(dz)f,h\rangle=\\
=\langle  Px(\varepsilon_m,\varepsilon_n)f,h\rangle,\quad f,h\in L_2,
\endgathered\eqno(1.11)
$$
where $P=v(N)\big(I+iv_1(N)\big)$ is a bounded normal operator that commutes with
$x(\varepsilon_m,\varepsilon_n)$. Put in (1.11) $f=\nu(t)$, $h=\nu(s)$, use
the generalized Schwarz inequality as above, and obtain, by (1.4), that
$$
\gathered
|\varPhi_{\varepsilon_m}(s,t)-\varPhi_{\varepsilon_n}(s,t)|=
|\langle \bsy\phi(\varepsilon_m,\varepsilon_n)(\nu(t)),\nu(s)\rangle|=\\
=|\langle x(\varepsilon_m,\varepsilon_n)(\nu(t)), P^*\nu(s)\rangle|\le\\
\le\langle  x(\varepsilon_m,\varepsilon_n)(\nu(t)),\nu(t)\rangle
\langle  x(\varepsilon_m,\varepsilon_n)PP^*(\nu(s)),\nu(s)\rangle\le\\
\le\|PP^*\|\langle x(\varepsilon_m,\varepsilon_n)(\nu(t)),\nu(t))\rangle
\langle x(\varepsilon_m,\varepsilon_n)(\nu(s)),\nu(s)\rangle=\\
=\|PP^*\|\left(X_{\varepsilon_m}(s,s)-X_{\varepsilon_n}(s,s)\right)
\left(X_{\varepsilon_m}(t,t)-X_{\varepsilon_n}(t,t)\right)
\endgathered\eqno(1.12)
$$
for all $(s,t)\in\R^2$ (in the next to the last inequality of (1.12), the
Reid inequality \cite[p. 59]{20} was used). From (1.12), by (1.10), it follows
that there exists $\varPhi'\in C(\R^2,\mathbb C)$ such that
$$
\lim_{n\to\infty}
\|\varPhi_{\varepsilon_n}-\varPhi'\|_{C(\R^2,\mathbb C)}=0.
$$
Now, exactly as for $X(s,t)$ and $X'(s,t)$, we can conclude in view of (1.6) and
(1.6)$^\prime$ that, for each $s\in\R$,
$$
\varPhi(s,t)=\varPhi'(s,t)
$$
for alwost every $t\in\R$, and conversely.
It follows that the function $\varPhi'(s,t)$ is a $K^0$-kernel of the integral
operator $\bsy\phi(N)$.
\end{proof}

\section{Representation of $K^0$-kernels by bilinear series}

\begin{theorem}
Let $N:D_N\to L_2$ be a normal integral operator induced by a
$K^0$-kernel $N(s,t)$ and having diagonal form
$$
Nf=\Sum_{n=1}^\infty\alpha_n\langle f,\varphi_n\rangle\varphi_n,\quad f\in D_N,
$$
where $\{\varphi_n\}_{n=1}^\infty\subset L_2$ is an orthonormal set,
the numbers $\alpha_n$ ($n=1,2,3,\dots$) are in a sector of angle less than
$\pi$ with vertex at $0$, and the series converges in $L_2$-norm to $Nf$.
Then the bilinear series
$$
\Sum_{n=1}^\infty\alpha_n\varphi_n(s)\adj{\varphi_n(t)}\eqno(2.0)
$$
converges both absolutely and uniformly in $\R^2$ to the kernel $N(s,t)$.
\end{theorem}
\begin{proof}
Take an $\alpha\in [0,2\pi]$ so that the eigenvalues
$\lambda_n=e^{i\alpha}\alpha_n=x_n+iy_n$ of the operator $P=e^{i\alpha}N$ 
are in the sector bounded by the rays $y=\pm lx$, where $x\ge0$, $l>0$.
The diagonal operator
$$
Tf=\dfrac{P+P^*}2f=\Sum_{n=1}^\infty x_n\langle f,\varphi_n\rangle\varphi_n,
\quad f\in D_T=D_N,\eqno(2.1)
$$
is an integral one having the $K^0$-kernel
$$
K(s,t)=\adj{K(s,t)}=\dfrac{e^{i\alpha}N(s,t)+e^{-i\alpha}\adj{N(t,s)}}{2}
$$
and the Carleman function $k(s)=\overline{K(s,\cdot)}$.
Assume, with no loss of generality, that $x_n>0$, $n=1,2,3\dots$.
Since $k\in C(\R,L_2)$, the equivalence class $Tf$ ($f\in D_T$) contains a
unique continuous function; this function has the form $\langle f,k(s)\rangle$ ($s\in\R$)
and belongs to $C(\R,\mathbb C)$. Therefore we can assume in what follows that
$$
\varphi_n\in C(\R,\mathbb C),\quad n=1,2,3\dots,\eqno(2.2)
$$
because $\dfrac1{x_n}T\varphi_n=\varphi_n$, $n=1,2,3\dots$.

Consider the functions
$$
K_m(s,t)=K(s,t)-\Sum_{n=1}^m x_n\varphi_n(s)\adj{\varphi_n(t)},\quad
m=1,2,3,\dots.
$$
By virtue of (2.2), these functions belong to $C(\R^2,\mathbb C)$.
Moreover, for every $f\in D_T$,
$$
\gathered
\Int_{\R}\Int_{\R}K_m(s,t)f(t)\adj{f(s)}\,dt\,ds
=\Sum_{n=1}^\infty x_n|\langle f,\varphi_n\rangle|^2-
\Sum_{n=1}^m x_n|\langle f,\varphi_n\rangle|^2=\\
=\Sum_{n=m+1}^\infty x_n|\langle f,\varphi_n\rangle|^2\ge 0.
\endgathered
$$
Therefore, $K_m(s,s)\ge 0$ for all $m$ and for all $s\in\R$ (see, for example,
\cite[p. 263]{17}), and hence
$$
K(s,s)\ge\Sum_{n=1}^\infty x_n|\varphi_n(s)|^2\quad\text{for all\ } s\in\R.
\eqno(2.3)
$$
By virtue of the Cauchy inequality,
$$
\left(\Sum_{n=p}^q x_n|\varphi_n(s)|\,|\varphi_n(t)|\right)^2\le
M\Sum_{n=p}^q x_n|\varphi_n(s)|^2\quad\text{for all\ } p, q\in\mathbb N,
\eqno(2.4)
$$
where $M=\max\limits_{s\in\R}K(s,s)$. Consequently, the series
$$
\Sum_{n=1}^\infty x_n\varphi_n(s)\adj{\varphi_n(t)}\eqno(2.5)
$$
absolutely converges on $\R^2$ and uniformly converges with respect to $t$.
Its sum-function $B(s,t)$ is continuous with respect to $t$ 
for each fixed $s$, and vice versa. Apply the Bessel inequality to the function
$k(s)\in L_2$ and obtain, for each fixed $s\in\R$,
$$
\Sum_{n=1}^\infty x_n^2|\varphi_n(s)|^2=
\Sum_{n=1}^\infty|\langle\varphi_n,k(s)\rangle|^2
\le\|k\|^2_{C(\R,L_2)}.
\eqno(2.6)
$$
It follows from the Riesz-Fisher theorem that $B(s,\cdot)\in L_2$ and that
$$
\Int_{\R}B(s,t)f(t)\,dt=\Sum_{n=1}^\infty x_n\langle f,\varphi_n\rangle\varphi_n(s)
\quad (f\in L_2)
$$
for all $s\in\R$.
The series on right-hand side of the last equation uniformly converges, since, by 
the Cauchy inequality and by (2.6),
$$
\gathered
\left|\Sum_{n=p}^q x_n\langle f,\varphi_n\rangle\varphi_n(s)\right|^2\le
\Sum_{n=p}^q x_n^2|\varphi_n(s)|^2\cdot\Sum_{n=p}^q|\langle f,\varphi_n\rangle|^2\le\\
\le\|k\|^2_{C(\R,L_2)}\Sum_{n=p}^q|\langle f,\varphi_n\rangle|^2.
\endgathered
$$
Note also that if $f\in D_T$, then by (2.1) the pointwise sum of this
series is the exactly continuous function $(Tf)(s)$. Therefore, one can write
$$
\Int_{\R}(K(s,t)-B(s,t))f(t)\,dt\equiv0\quad\text{on $\R$}
$$
for each $f$ from the dense set $D_T$ in $L_2$.
Consequently, for every fixed $s\in\R$,
$$
K(s,t)=B(s,t)\eqno(2.7)
$$
for all $t\in\R$, since both functions are continuous with respect to $t$.
In particular, for each $s\in\R$,
$$
K(s,s)=B(s,s)=\Sum_{n=1}^\infty x_n|\varphi_n(s)|^2.\eqno(2.8)
$$
The series (2.8) can be considered on the compactum $[0,\infty]$, keeping in
mind that $K(\infty,\infty)=\varphi_n(\infty)=0$, $n=1,2,3,\dots$. Hence, by
the Dini theorem the series (2.8) converges uniformly on $\R$. From (2.4) it follows
that the series (2.5) uniformly converges on $\R^2$ to the kernel $K(s,t)$. Using
the estimate $|\alpha_n|=|\lambda_n|\le x_n\sqrt{1+l^2}$ $(n=1,2,3,\dots)$
and the Cauchy inequality, we have
$$
\Sum_{n=p}^q|\lambda_n||\varphi_n(s)||\varphi_n(t)|\le\sqrt{1+l^2}
\left(\Sum_{n=p}^q x_n|\varphi_n(s)|^2\right)^{1/2}\left(\Sum_{n=p}^q 
x_n|\varphi_n(t)|^2\right)^{1/2}.
$$
Therefore, the series (2.0) uniformly converges on $\R^2$. Its sum $\Omega(s,t)$
is a function belonging to $С(\R^2,\mathbb C)$. Now apply preceding reasoning
to the functions $\Omega(s,t)$ and $N(s,t)$ in place of $B(s,t)$ and $K(s,t)$
in (2.7) to infer that $\Omega(s,t)\equiv N(s,t)$ on $\R^2$.
\end{proof}

\bibliographystyle{amsplain}

\end{document}